\documentclass[aap,seceqn,dvips]{arximspdf}

\doi{10.1214/105051606000000880}
\volume{17}
\issue{2}
\pubyear{2007}
\firstpage{716}
\lastpage{740}

\makeatletter
\newproclaim{res}{Result}[section]
\newtheorem{theorem}[res]{Theorem}
\newproclaim{rem}[res]{Remark}
\newtheorem{prop}[res]{Proposition}
\newtheorem{lem}[res]{Lemma}
\newtheorem{cor}[res]{Corollary}
\newproclaim{defi}[res]{Definition}

\newproclaim{assum}[res]{Assumption}
\let\widebar\overline
\makeatother

\begin{document}
\begin{frontmatter}

\title{On the density of properly maximal claims in financial markets
with transaction costs\thanksref{T1}}
\runtitle{On the Density of Maximal Claims}
\thankstext{T1}{Supported by the grant ``Distributed Risk
Management'' in the Quantitative Finance initiative funded by EPSRC and
the Institute and Faculty of Actuaries.}

\begin{aug}
\author[A]{\fnms{Saul} \snm{Jacka}\corref{}\ead[label=e1]{s.d.jacka@warwick.ac.uk}} and
\author[A]{\fnms{Abdelkarem} \snm{Berkaoui}\ead[label=e2]{a-k.berkaoui@warwick.ac.uk}}
\runauthor{S. Jacka and A. Berkaoui}
\affiliation{University of Warwick}
\address[A]{Department of Statistics\\
University of Warwick\\
Coventry CV4 7AL\\
UK\\
\printead{e1}\\
\phantom{\textsc{E-mail}:} \printead*{e2}} 
\end{aug}

\received{\smonth{3} \syear{2006}}
\revised{\smonth{10} \syear{2006}}

%
\begin{abstract}
We consider trading in a financial market with proportional transaction
costs. In the frictionless case, claims are maximal if and only if they
are priced by a consistent price process---the equivalent of an
equivalent martingale
measure. This result fails in the presence of transaction costs.
A properly maximal claim is one which does have this property.
We show that the properly maximal claims are dense in the set of
maximal claims (with the topology of convergence in
probability).\looseness=-1
\end{abstract}

%

\end{frontmatter}

\section{Introduction}\label{s1}
We consider a discrete-time market in $d$ assets with transaction
costs. We suppose that $\mathcal A$ is the cone of claims attainable
from $0$ by trading. In~\cite{JBW}, following on from Schachermayer
\cite{schacher}, Kabanov \cite{kaban}, Kabanov, Stricker and Rasonyi
\cite{kaban2} and \cite{kaban3} and many others, Jacka, Berkaoui and
Warren showed that if $\mathcal A$ is arbitrage-free (i.e., contains no
positive elements) then whilst $\mathcal A$ may not be closed, its
closure [in $L^0(\mathbb R^d)$] is also a cone of attainable claims
under a new price system and is arbitrage-free if and only if there is
a consistent price process for $\mathcal A$ (Theorems 3.6 and 4.12).
Here, a consistent price process is essentially given by a strictly
positive element in the polar cone of $\mathcal A\cap L^1$. A
consistent price process is a suitable generalization of the concept of
the density of an equivalent martingale measure (EMM).\looseness=-1

Given a claim $X\in\mathcal A$, a standard question is how to hedge it.
In other words, how to find a self-financing trading strategy which
achieves a final portfolio of $X$ with 0 initial endowment. In the
context of frictionless trading, this is achieved by seeking maximal
claims---claims $Y$ which are maximal in $\mathcal A$ with respect to
the partial order
\[
W\geq X\Leftrightarrow W-X\in\mathbb R^d_+\qquad\mbox{a.s.},
\]
(see \cite{DS1,DS2,DS3}). It follows from Kramkov's celebrated result
on optional decompositions (\cite{kram}) that, at least in a
discrete-time context, a claim $X$ in $\mathcal A$ is maximal if and
only if it is priced at 0 by some EMM. It also follows that this is
true if and only if $[\mathcal A,X]$, the cone generated by $\mathcal
A$ and $-X$, is arbitrage-free, in which case its closure is also
arbitrage-free.

Consequently (see \cite{jacka} or \cite{DS1}), one may obtain a hedging strategy for a
maximal claim by martingale representation.

Regrettably, when there are transaction costs,
just as $\mathcal A$ may be arbitrage-free but $\bar\mathcal A$
contain an arbitrage, so, in this context,
a claim $X$ may be maximal and yet the
closure of $[\mathcal A,X]$ contain an arbitrage.

In the language of optimization theory, a maximal claim such that the
closure of $[\mathcal A,X]$ is arbitrage-free is said to be proper
efficient with respect to $L^{0,+}$. We shall refer to such claims as
properly maximal. We shall show in Theorem \ref{equiv} that a properly
maximal claim is priced by  some consistent price process and that
martingale representation can be used to obtain a hedging strategy. It
is then of interest (for hedging purposes) as to whether one can
approximate maximal claims by properly maximal claims. This is a
problem with a long and distinguished history in optimization theory,
going back to~\cite{ABB}. We give a positive answer (up to
randomization) in Theorem~\ref{t1manu}: the collection of properly
maximal claims is dense in the set of maximal claims.\looseness=-1

In a continuous time framework, the problem is more delicate. Indeed,
the task of defining a notion of admissible trading strategy, that has
a meaningful financial interpretation, is still in progress. A first
solution has been given by Kabanov \cite{kaban}, Kabanov and Last
\cite{KabLast} and Kabanov and Stricker \cite{KabSt}, where the
efficient friction assumption was made. More precisely, an admissible
self-financing trading strategy was defined as an adapted,
vector-valued, c\'adl\'ag process of finite variation whose increments
lie in the corresponding trading/solvency cones and whose terminal
value is bounded from below by a constant with respect to the order
induced by the terminal solvency cone. Campi and Schachermayer
\cite{CScha} extend these results to bid-ask processes which are not
necessarily continuous. In this framework, the discrete-time
methodology cannot be adopted, as it is based on the fact that the cone
of attainable claims for zero endowment is a \textit{finite} sum of
one-period trading cones.
\vfill\eject

\section{Background, notation and preliminary results}\label{back}

\subsection{Efficient and proper efficient points}
Given a topological vector space $Z$, a pointed, closed, convex cone
$C$ defines a partial order $\buildrel C\over\leq $ on $Z$ by
\[
x\buildrel C\over\leq y\Leftrightarrow y-x\in C.
\]

For a subset $B\subset Z$, we denote by $\operatorname{cone}(B)$ the
cone generated by $B$, that is,
\[
\operatorname{cone}(B)=\{\lambda b\dvtx  \lambda\in\mathbb R^+, b\in
B\}.
\]
For a convex set $D\subset Z$, we denote by $\operatorname{lin}(D)$ the
lineality subspace of $D$:
\[
\operatorname{lin}(D)= \bigcup_{\mathrm{subspaces}\ V\subseteq D}V,
\]
and recall that if $D$ is also a cone then
$\operatorname{lin}(D)=D\cap(-D)$.

\begin{defi}
Given a subset
$\mathcal A\subset Z$, we say that
\[
\theta \in\mathcal A\mbox{ is $C$-efficient if }
\operatorname{cone}(\mathcal A-\theta )\cap C=\{0\},
\]
and
\[
\theta \in\mathcal A\mbox{ is proper $C$-efficient if }
\overline{\operatorname{cone}(\mathcal A-\theta )}\cap C=\{0\}.
\]
If the cone $C$ is not pointed then we change the definitions as
follows:
\[
\theta \in\mathcal A\mbox{ is $C$-efficient if }
\operatorname{cone}(\mathcal A-\theta )\cap C\subset
\operatorname{lin}(C),
\]
and
\[
\theta \in\mathcal A\mbox{ is proper $C$-efficient if }
\overline{\operatorname{cone}(\mathcal A-\theta )}\cap C\subset
\operatorname{lin}(C).
\]
\end{defi}

One of the main problems in multi-criteria optimization theory is to
show that each efficient point can be approximated by a sequence of
proper efficient points---hereafter we refer to this as the
\textit{density problem}.

\begin{rem}
It is easy to see that if $\theta \in\mathcal A$ is $C'$-efficient,
with $C'$ a pointed closed convex cone such that
$C\setminus\{0\}\subset \operatorname{int}(C')$, then $\theta $ is also
proper $C$-efficient.

Take $X\in\overline{\operatorname{cone}(\mathcal A-\theta )}\cap C$. If
$X\neq0$ then $X\in \operatorname{int}(C')$ so we can take a
neighborhood of $X$, $U$, such that $0\notin U$, $U\cap
\operatorname{cone}(\mathcal A-\theta)\neq\varnothing$ and $U\subset
C'$. But this implies that there is a $y$ with $y\neq0$, $y\in U\cap
\operatorname{cone}(\mathcal A-\theta)$ and $y\in C'$ which contradicts
the $C'$-efficiency of $\theta$.

The inverse implication is not always true unless we suppose further
conditions on the triplet $(Z,C,\mathcal A)$.
\end{rem}

\begin{rem}
One way then to solve the density problem is to construct a sequence of
pointed closed convex cones $(C_n)_{n\geq1}$ which decrease to the
convex cone $C$ and are such that $C\setminus\{0\}\subset
\operatorname{int}(C_n)$. Such a sequence is called a
\textit{$C$-approximating sequence} (or family) of cones. In this case
the set $(\theta +C_n)\cap\mathcal A$ will converge to the set $(\theta
+C)\cap\mathcal A$ which is reduced to the singleton $\{\theta \}$ if
$\theta $ is $C$-efficient. In consequence, each $C_n$-efficient point
$\theta _n$ is proper $C$-efficient and any sequence of $C_n$-efficient
points will converge to $\theta $.
\end{rem}

In \cite{ABB}, Arrow, Barankin and Blackwell solved the density problem
in the finite-dimensional case where: $Z=\mathbb R^n$, $C=\mathbb
R^n_+$ with $\mathcal A$ a compact, convex set in~$\mathbb R^n$. This
theorem was extended to cover more general topological vector spaces
(see \cite{BM,Borw-Zhang,Gong,Hartley,Henig,Ster-K,wan}). In
\cite{Ster-K}, Sterna-Karwat proved that in a normed vector space~$Z$,
there exists a $C$-approximating sequence of cones if and only if
\[
C^{+, i}\buildrel \mathrm{def}\over=\bigl\{\lambda \in Z^*; \lambda >0
\mbox{ on } C\setminus\{0\} \bigr\}\neq\varnothing.
\]
She then applied this result to the density problem for a compact
convex set $\mathcal A$.

The case of a locally convex vector space was discussed by, among
others, Fu Wantao in \cite{wan}. He solved the density problem by
supposing that the convex cone $C$ admits a \textit{base $B$}. This
means that $B$ is a convex set, $0\notin \widebar{B}$ and
$C=\overline{\operatorname{cone}(B)}$. He used this assumption to
construct a $C$-approximating family of cones.

For a more recent survey of such techniques see \cite{dan}.

\subsection{Notation and further background}
We are equipped with a filtered probability space $(\Omega,{\mathcal
F},({\mathcal F}_t)_{t=0,\ldots,T},{\mathbb P})$. We denote the
real-valued \mbox{${\mathcal F}_t$-}mea\-sur\-able functions by $m{\mathcal
F}_t$, the nonnegative subset by $m{\mathcal F}_t^+$, the bounded
real-valued \mbox{${\mathcal F}_t$-}measurable functions by $b{\mathcal
F}_t$ and the nonnegative subset by $b{\mathcal F}_t^+$. We denote the
space of ${\mathcal F}$-measurable random variables in $\mathbb R^d$ by
$L^0({\mathcal F},\mathbb R^d)$ (with the metric which corresponds to
the topology of convergence in measure) or just $L^0_d$. And we denote
the almost surely nonnegative and nonpositive subsets by $L^{0,+}_d$
and $L^{0,-}_d$ respectively.

\begin{defi}
Through the paper we adopt the following notation. Suppose $F\in
{\mathcal F}$,
$D\subset L^0_d$ is a convex cone and $\xi\in D$. We define
\[
[D,\xi]\buildrel \mathrm{def}\over=\operatorname{cone}(D-\xi).
\]
Notice that, because $D$ is a cone this satisfies
\[
[D,\xi]=D-\mathbb R^+\xi.
\]
Moreover, $[D,\xi]$ inherits the convexity property from $D$.

We define
\[
\xi(F)\buildrel \mathrm{def}\over=\xi{\mathbf1}_F,
\]
where ${\mathbf 1}_F$ is the indicator function of $F$, and
\[
D(F)\buildrel \mathrm{def}\over=\{\xi(F)\dvtx  \xi\in D\}.
\]
We say that $D$ is \textit{arbitrage-free} if $D\cap L^{0,+}_d=\{0\}$.
We denote the complement of a subset $B$ in $\Omega $ by $B^c$.
\end{defi}

We recall the setup from Schachermayer's paper \cite{schacher}: we may
trade in $d$ assets at times $0,\ldots,T$. We may burn any asset and
otherwise trades are on terms given by a bid-ask process $\pi $ taking
values in $\mathbb R^{d\times d}$, with $\pi $ adapted. The bid-ask
process gives the (time $t$) price for one unit of each asset in terms
of each other asset, so that
\[
\pi ^{i,i}_t=1\qquad \forall i,
\]
and $\pi ^{i,j}_t$ is the (random) number of units of asset $i$ which
can be traded for one unit of asset $j$ at time $t$. We assume (with
Schachermayer) that we have  ``netted out'' any advantageous trading
opportunities, so that, for any $t$ and any $i_0,\ldots,i_n$:
\[
\pi^{i_0,i_n}_t\leq\pi^{i_0,i_1}_t\cdots\pi^{i_{n-1},i_n}_t.
\]

The time $t$ trading cone, $K_t$, consists of all those random
trades (including the burning of assets) which are available at
time $t$. Note that Schachermayer refers to this cone as $-\hat K_t$.
Thus we can think of $K_t$ as consisting of all those
random vectors which live (almost surely) in a random closed convex
cone $K_t(\omega)$.

Denoting the $i$th canonical basis vector of ${\mathbb R}^d$ by $e_i$,
$K_t(\omega)$ is defined as the finitely-generated (hence closed)
convex cone with generators $\{e_j-\pi ^{i,j}_t(\omega)e_i,1\leq i\neq
j\leq d;\mbox{ and }-e_k,1\leq k\leq d\}$. The reader is referred to
Theorem 4.5 and the subsequent Remark 4.6 of \cite{JBW}.

We shall say that $\eta$ is a \textit{self-financing process} if
$\eta_t-\eta_{t-1}\in K_t$ for each $t$, with $\eta_{-1}\buildrel
\mathrm{def}\over=0$. We say that $\underline\xi$ is a \textit{hedging
strategy} if $\underline\xi\in K_0\times\cdots\times K_T$.

It follows that the cone of claims attainable from zero endowment is
$K_0+\cdots+K_T$ and we denote this by $\mathcal A$. As we said in the
\hyperref[s1]{Introduction}, $\mathcal A$ may be arbitrage-free and yet its closure
may contain an arbitrage. However, by Theorem 3.6 of~\cite{JBW}, we may
and shall assume that if $\bar\mathcal A$ is arbitrage-free then (by
adjusting the bid-ask process) $\mathcal A$ is closed and
arbitrage-free. We should remark at this point that a very small
generalization of this theorem allows us to continue to make this
assumption merely if each $K_t$ is a finitely ${\mathcal
F}_t$-generated convex cone with the ${\mathcal F}_t$-measurable
generators given by $\Pi^1_t,\ldots,\Pi^n_t$: that is,
\[
K_t = \Biggl\{\sum_{i=1}^n\alpha_i\Pi^i_t\dvtx  \alpha_i\in m{\mathcal
F}_t^+ \Biggr\}.
\]
Henceforth, any such cone will be described as a \textit{finitely
${\mathcal F}_t$-generated cone}.

For any decomposition of $\mathcal A$ as a sum of convex cones:
\[
\mathcal A=M_0+\cdots+M_t ,
\]
we call elements of $M_0\times\cdots\times M_t$, which almost surely
sum to 0, \textit{null-strategies} (with respect to the decomposition
$M_0+\cdots+M_t$). We denote the set of null-strategies by $\mathcal
N(M_0\times\cdots\times M_t)$. For convenience we denote
$K_0\times\cdots\times K_T$ by~${\mathcal K}$.

In what follows we shall often use (a slight generalization of)
Schachermayer's key result (Remark 2.8 after the proof of Theorem 2.1
of \cite{schacher}):

\begin{lem}\label{s}
Suppose that
\[
\mathcal A=M_0+\cdots+M_{t-1}+M_{t}
\]
is a decomposition of $\mathcal A$ into convex cones with $M_s\subseteq
L^0({\mathcal F}_s,\mathbb R^d)$ for $0\leq s\leq t-1$, and $b{\mathcal
F}_s^+ M_s\subseteq M_s$ for each $s\leq t$. If $\mathcal
N(M_0\times\cdots\times M_t)$ is a vector space and each $M_t$ is
closed, then $\mathcal A$ is closed.
\end{lem}

\begin{rem}
Theorem 3.6 of \cite{JBW} establishes that where $\bar\mathcal A$ is
arbitrage-free, the revised bid-ask process gives rise to finitely
generated cones $\tilde K _t$ with the further property that $\mathcal
N(\tilde{\mathcal K})$ is a vector space.
\end{rem}

\begin{cor}
\label{close} There exists a family $(M_t)_{t=0,\ldots,T}\subset
L^0_d$ with each $M_t$ a closed convex cone ${\mathcal F}_t$-generated
by a
finite family of $\mathbb R^d$-valued ${\mathcal F}_t$-measurable
vectors, such that:
\[
\mathcal A=M_0+\cdots+M_T
\]
and
\[
\mathcal N(M_0\times\cdots\times M_T )=\{\underline{0}\}.
\]
\end{cor}

\begin{pf}
We have assumed that
\[
\mathcal A=\tilde K_0+\cdots+\tilde K_T,
\]
and that $\mathcal N(\tilde{\mathcal K})$ is a vector space. This
implies (by Lemma \ref{s}) that $\mathcal A$ is closed and that $
\eta_t\buildrel \mathrm{def}\over=\mathcal N (\tilde
K_t\times\cdots\times\tilde K_T )$ is a vector space for each
$t=0,\ldots,T-1$. Define $\rho_t$ to be the projection of the closed
vector space $\eta_t$ onto its first component and define $M_t=\tilde
K_t\cap\rho_t^\perp$ with $M_T=\tilde K_T$. We verify easily that the
family $(M_0,\ldots,M_T)$ satisfies the conditions of the corollary.
\end{pf}

So, from now on we make the following:

\begin{assum}\label{neat}
The cone of claims attainable from 0, $\mathcal A$, can be written as
$\mathcal A={K_0+\cdots+K_T}$ where each $K_t$ is a finitely ${\mathcal
F}_t$-generated convex cone and $\mathcal N({\mathcal
K})=\{\underline{0}\}$. Consequently $\mathcal A$ is closed.
\end{assum}

In what follows, the terms  ``maximality'' and  ``proper maximality''
are defined with respect to the cone $C=L^{0,+}_d$. For more general
ordering cones we continue to use the terms  ``efficiency'' and
``proper efficiency.''

\subsection{Maximal claims and representation}
Recall that in the frictionless setup, $X\in\mathcal A$ is maximal if
and only if there is an EMM ${\mathbb Q}$ such that $\mathbf E_{\mathbb
Q}X=0$. Moreover, in that case, denoting the collection of EMMs by $Q$,
\[
\mathbf E_{\mathbb Q}X=0
\]
for every ${\mathbb Q}\in Q$. In this case, defining $V_t$ as the
common value of $\mathbf E_{\mathbb Q}[X|{\mathcal F}_t]$, the process
$V$ is a $Q$-uniform martingale and hence is representable as a
stochastic integral with respect to the discounted price process. See
\cite{DS1}~and~\cite{jacka} for details. The stochastic integrand
essentially then gives a hedging strategy.

Recall from \cite{schacher} that the concept of consistent price
process is a suitable generalization of the concept of EMM. To be
precise, a consistent price process is a martingale $(Z_t)_{0\leq t\leq
T}$, with $Z_t$ taking values in $K^*_t\setminus\{0\}$, where $K^*_t$
is the polar cone (in $\mathbb R^d$) of $K_t$. The value $Z_t$ plays
the same role as the density of the restriction of an EMM to ${\mathcal
F}_t$ in the frictionless setting.

\begin{theorem}\label{equiv}
In the case of transaction costs, suppose that $X\in\mathcal A$, then:
\begin{longlist}[(2)]
\item[(1)] there exists a consistent price process, $Z$, such that
$\mathbf EZ_T\cdot X=0$ if and only if $X$ is properly maximal.

\item[(2)] Suppose that $X$ is properly maximal and let $Q^Z$ be the
collection of EMMs for the consistent price process $Z$. Then $V^Z$,
defined by
\[
V^Z_t=\mathbf E_{\mathbb Q}[Z_T\cdot X|{\mathcal F}_t],
\]
is independent of the choice of ${\mathbb Q}\in Q^Z$ and is a
$Q^Z$-uniform martingale.
\end{longlist}
\end{theorem}

We recall first Theorem 4.12 of~\cite{JBW}, that we will need in the
next proof:
\begin{itemize}
\item[] \textit{$\bar\mathcal A$}, \textit{the closure of $\mathcal A$
in $L^0$}, \textit{is arbitrage-free iff there is a consistent price
process~$Z$}. \textit{In this case}, \textit{for every strictly
positive ${\mathcal F}_T$-measurable $\phi\dvtx \Omega
\rightarrow(0,1]$ we may find a consistent price process $Z$ such that
$|Z_T|\leq c\phi$ for some positive constant~$c$}.
\end{itemize}

\begin{pf*}{Proof of Theorem~\ref{equiv}}
First recall Theorem 4.16 of \cite{JBW}:

If $\theta\in L^0({\mathcal F}_T,\mathbb R^d)$ and $\mathcal A$ is
closed and arbitrage-free, the following are equivalent:
\begin{longlist}[(ii)]
\item[(i)] $\theta\in\mathcal A$.

\item[(ii)]For every consistent pricing process $Z$ such that the
negative part $(\theta\cdot Z_T)^-$ of the random variable $\theta\cdot
Z_T$ is integrable, we have
\[
\mathbf E[\theta\cdot Z_T]\leq0.
\]
\end{longlist}

Proof of (1): Suppose $X\in\mathcal A$ and $Z$ is a consistent price
process and $X$ and $Z$ satisfy condition (ii) above. Write $X$ as
%
\begin{equation}\label{xdec}
X=\sum _{s=0}^T\xi_s,
\end{equation}
with $\xi_s\in K_s$ for each $s$, and then, for each $t$, denote $\sum
_{s=0}^t\xi_s$ by $X_t$.

Notice that, since $Z$ is consistent, $Z_t$ is in $K_t^*$ and so
$Z_t\cdot \xi_t\leq0$ for each $t$. So, in particular,
\[
Z_T\cdot X=Z_T\cdot (X_{T-1}+\xi_T)\leq Z_T\cdot X_{T-1}
\]
and so
\[
(Z_T\cdot X)^-\geq(Z_T\cdot X_{T-1})^-,
\]
and thus $(Z_T\cdot X_{T-1})^-$ is integrable. Now, for each $t$,
$X_t\in \mathcal A$ (since it is in $K_0+\cdots+K_t$) so, by part (ii)
of Theorem 4.16 of \cite{JBW},
\[
\mathbf E[Z_T\cdot X_{T-1}]\leq0
\]
and so, in particular, $Z_T\cdot X_{T-1}$ is integrable and
\[
\mathbf E[Z_T\cdot X]\leq\mathbf E[Z_T\cdot X_{T-1}]=\mathbf
E\bigl[\mathbf E[ Z_T\cdot X_{T-1}|{\mathcal F}_{T-1}]\bigr]=\mathbf E[Z_{T-
1}\cdot X_{T-1}].
\]
Now we iterate the argument [which we may do since $\mathcal
A_t\buildrel \mathrm{def}\over=K_0+\cdots+K_t$ is closed for each $t$,
which follows from our Assumption \ref{neat} and $(Z_0,\ldots,Z_t)$ is
a consistent price process for $(K_0\times\cdots\times K_t)$]. We see
that if $X$ and $Z$ satisfy the conditions of (ii) above then
%
\begin{equation}\label{rep2}
\mathbf E[Z_T\cdot X]=\sum _{s=0}^T\mathbf E[Z_s\cdot\xi_s].
\end{equation}

Now, suppose that $Z$ is a consistent price process for $\mathcal A$,
$X\in\mathcal A$ and\break \mbox{$\mathbf EZ_T\cdot X=0$}. Recalling equation
(\ref{xdec}), it follows from the consistency of $Z$ and (\ref{rep2})
that
\[
Z_t\cdot\xi_t= 0,
\]
for each $t$. Now it is easy to check that $\overline{[\mathcal
A,X]}=\overline{K_0^{\xi_0}+\cdots+K_T^{\xi _T}}$, where $K_t^x$
denotes the ${\mathcal F}_t$-cone obtained from $K_t$ by adding the
generator $-x$. We apply Theorem 4.12 of~\cite{JBW}, with $\bar
\mathcal A$ replaced by $\overline{[\mathcal A,X]}$ to conclude that
$\overline{[\mathcal A,X]}$ is arbitrage-free.

Conversely, suppose that $\overline{[\mathcal A,X]}$ is arbitrage-free
then, by Theorem 4.12 of \cite{JBW} again, there is a consistent price
process, $Z$, for $K_0^{\xi_0}+\cdots+K_T^{\xi_T}$ satisfying
\[
|Z_T|\leq\frac{c}{(1+\sum _{s=0}^T|\xi_s|_{\mathbb R^d})}
\]
for some positive $c$. Notice that, since $Z_t\in(K_t^{\xi_t})^*$ and
both $\xi_t$ and $-\xi_t$ are in $K_t^{\xi_t}$, we must have
$Z_t\cdot\xi_t=0$ for each $t$. It follows, a fortiori, from
the fact that $Z$ is consistent for $K_0^{\xi_0}+\cdots+K_T^{\xi_T}$
that $Z$ is consistent for $K_0+\cdots+K_T$. Notice that
$Z_t\cdot\xi_t$ is bounded by $c$, and hence integrable, for each $t$
and so, by the usual arguments
\[
\mathbf EZ_T\cdot X=\sum _{t=0}^T \mathbf E[Z_t\cdot\xi_t]=0.
\]

To prove (2), simply notice that $V^Z_t=\mathbf E_{\mathbb Q}[Z_T\cdot
X|{\mathcal F}_t]=Z_t\cdot X_t$ by virtue of the usual tower-property
arguments and the fact that ${\mathbb Q}$ is an EMM for $Z$.
\end{pf*}

\begin{rem}
Of course, representation does not guarantee that the  ``hedging
strategy'' $\xi$ is admissible: it may be that it is  ``priced at 0''
by $Z$ but is still not in $\mathcal A$ because some other consistent
price process assigns it a positive price. This can happen if $\xi_t$
is in $\operatorname{span}(K^{\xi_t})$ but not in $K^{\xi_t}$.
\end{rem}

\section{An example of a maximal claim which is not proper}\label{s3}
We take a simple setup for trading in two assets over two time periods.
We set $T=1$, $d=2$ and $\Omega=\mathbb N$. We take ${\mathcal F}_0$ as
the trivial $\sigma$-algebra, set ${\mathcal F}_1=2^\mathbb N$ and
define ${\mathbb P}$ as any probability measure which puts positive
mass on each point of $\Omega$. The bid-ask process $\pi$ is given by:
\[
\pi^{1,2}_0=1,\qquad \pi^{2,1}_0=k,\qquad
\pi^{2,1}_1=2\quad\mbox{and}\quad \pi^{1,2}_1=k,
\]
where $k$ will be taken suitably large.

The claim $\theta$ is defined by
\[
\theta\buildrel
\mathrm{def}\over=\biggl(1-\frac{1}{\omega}\biggr)e_2-\biggl(1-\frac
{1}{2\omega}\biggr)e_1,
\]
which corresponds to the following trading strategy---at
time 0 buy 1 unit of asset~2 for one unit of asset 1. At time 1 sell
$\frac{1}{\omega}$ units of asset 2 for $\frac{1}{2\omega}$ units of
asset 1.

\begin{prop}
If $k$ is sufficiently large, the claim $\theta$ defined above is
maximal but not properly maximal.
\end{prop}

\begin{pf}
Recall from  \cite{schacher} that a {\it strictly} consistent price
process is a martingale $(Z_t)_{0\leq t\leq T}$ with $Z_t$ taking
values in $\operatorname{rint}(K^*_t)\setminus\{0\}$, where $\operatorname{rint}$ denotes relative
interior. Theorems 1.7 and 2.1 of \cite{schacher} then show that if
there is a {\it strictly} consistent price process then $\mathcal A$ is
closed and arbitrage-free.

Now note first that if we take $Z_t=(Z^1,Z^2)=(1,\frac{3}{4})$ for both
$t=0$ and $t=1$ then $Z$ is a strictly consistent price process because
it is clear that $Z_t$ lies in the interior of $K_t^*$ for each $t$.
This follows since $0<\frac{Z^2}{Z^1}=\frac{3}{4}<\pi^{1,2}$ and
$0<\frac{Z^1}{Z^2}=\frac{4}{3}<\pi^{2,1}$.

Now to show that $\theta$ is maximal, suppose that $\phi\in\mathcal
A$ and
$\phi\geq\theta$ a.s. Let
\[
\phi=\xi_0+\xi_1,
\]
where $\underline\xi\equiv(\xi_0,\xi_1)\in K_0\times K_1$. It is clear
that we may suppose without loss of generality that $\xi_0$ is either
some positive multiple of $e_2-e_1$ or of $e_1-ke_2$. Similarly we may
suppose that $\xi_1$ is either some positive ${\mathcal
F}_1$-measurable multiple of $e_1-2e_2$ or of $e_2-ke_1$. By taking $k$
sufficiently large we may rule\ out the second possibility in each
case. This leaves us with the case where for suitable $a\in\mathbb R^+$
and $B\in m{\mathcal F}_1^+$
\begin{eqnarray*}
\xi_0&=& a(e_2-e_1),
\\
\xi_1&=& B(e_1-2e_2)
\end{eqnarray*}
and
\[
\phi=(a-2B)e_2-(a-B)e_1.
\]
Now if $\phi\geq\theta$ a.s. then we must have (comparing coefficients
of $e_2$ in $\phi$ and $\theta$):
%
\begin{equation}\label{ex1}
a-2B(n)\geq1-\frac{1}{n}\qquad\mbox{for all }n.
\end{equation}
Taking $\limsup_{n\rightarrow\infty}$ in (\ref{ex1}) we see that we must
have $a\geq1$.

However, comparing coefficients of $e_1$, we must have
%
\begin{equation}\label{ex2}
2\bigl(B(n)-a\bigr)\geq2\biggl(\frac{1}{2n}-1\biggr)\qquad\mbox{for all
}n,
\end{equation}
and adding (\ref{ex1}) and (\ref{ex2}), we see that we must
have $a\leq1$. Hence we see that $a=1$ and so, from (\ref{ex1}) and
(\ref{ex2}), $B(n)=\frac{1}{2n}$
and $\phi=\theta$. This establishes the maximality of $\theta$.

Now we shall show that $\mathcal A^\theta\buildrel
\mathrm{def}\over=\overline {\operatorname{cone}(\mathcal A-\theta
)}=\overline{[\mathcal A,\theta]}$ contains an arbitrage and hence that
$\theta$ is not proper.

Let $\xi$ denote the strategy above which attains $\theta$, so
\[
\xi_0=(e_2-e_1)
\]
and
\[
\xi_1=\frac{1}{2\omega}(e_1-2e_2).
\]
Notice that, since $\xi_0\in\mathcal A$, $\psi\buildrel
\mathrm{def}\over=-\xi_1=\xi _0-\theta\in\mathcal A^\theta$. It
follows, since $(\frac{1}{2}e_1-e_2)\in K_1$, that
$(\frac{n}{\omega}-1)\mathbf1_{(\omega\leq n)}(\frac
{1}{2}e_1-e_2)\in\mathcal A$ and so (adding $n\psi$)
\[
x_n\buildrel
\mathrm{def}\over=\bigl(e_2-\tfrac{1}{2}e_1\bigr)\mathbf1_{(\omega\leq
n)}\in\mathcal A^\theta.
\]
Letting $n\rightarrow\infty$ we deduce (from the closedness of
$\mathcal A^\theta$) that $(e_2-\frac{1}{2}e_1)\in\mathcal A^\theta
$ and so, adding
$e_1-e_2=-\xi_0=\xi_1-\theta\in\mathcal A^\theta$, we see that
$\frac{1}{2}e_1\in\mathcal A^\theta$, which is an arbitrage.
\end{pf}

\section{Some general results and the case $\operatorname{lin}(A)=\{0\}$}

\begin{defi}
We denote by $\mathcal A_{t,T}$ the closed cone $K_t+\cdots+K_T$. Given
$\theta\in\mathcal A$, we denote by $\theta_{t,T}$ the sum
$\theta_t+\cdots+\theta_T$.

We say that the decomposition of $\theta $:
\[
\theta=\theta_0+\cdots+\theta_T
\]
is a \textit{special decomposition} if, for each $t=0,\ldots, T-1$,
\[
\theta_t \mbox{ is efficient in } K_t\cap(\theta_{t,T}-\mathcal
A_{t+1,T})\mbox{ with respect to }-\!\mathcal A_{t+1,T}.
\]
To be more explicit, the decomposition is special if, for each
$t=0,\ldots, T-1$,
\begin{longlist}[(1)]
\item[(1)] $z\in\mathcal A_{t+1,T}$

\item[]and

\item[(2)] $\theta_t-z\in K_t$

\item[]together imply that

\item[(3)] $z=0$ if $\operatorname{lin}(\mathcal A_{t+1,T})=\{0\}$ or,
more generally, $z\in \operatorname{lin}(\mathcal A_{t+1,T})$.
\end{longlist}
\end{defi}

\begin{rem}
Notice that a special decomposition is in $\mathcal K$.
\end{rem}
\begin{rem}
If we think of hedging a claim, then a special decomposition is
``lazy'' in that it defers taking action until as late as possible, in
some sense.
\end{rem}

\begin{rem}
In the example of Section \ref{s3}, the decomposition
\[
(\theta_0,\theta_1)=\bigl(e_2-e_1,(1/2\omega)(e_1-2e_2)\bigr),
\]
is a special decomposition of $\theta=\theta_0+\theta_1$.

To show this, we want to prove that $\theta_0$ is efficient in
$K_0\cap(\theta-K_1)$ with respect to the order generated by $-K_1$. So
let $\xi_0\in K_0\cap(\theta-K_1)$ be such that
$\eta_1=\theta_0-\xi_0\in K_1$. Then there exists $a_0,b_0\in\mathbb
R^+$ and $a_1,b_1\in m{\mathcal F}^+_1$ such that
\[
\xi_0=a_0(e_2-e_1)+b_0(e_1-ke_2)
\]
and
\[
\eta_1=a_1(e_1-2e_2)+b_1(e_2-ke_1).
\]
We deduce that
\[
b_1=\dfrac{1-a_0+(2-k)b_0}{2k-1}
\]
and
\[
a_1=\dfrac{(k-1)(-1+a_0-(k+1)b_0)}{2k-1}.
\]
We take $k\geq10$ and since $a_1,b_1\geq0$ we obtain that $a_0=1$
and $b_0=0$, which means that $\eta_1=0$ and $\theta_0=\xi_0$. This
establishes the desired efficiency of $\theta_0$.
\end{rem}

We shall now show that every claim in $\mathcal A$ has a special decomposition.

\begin{theorem}
\label{decomp}
Given $\theta\in\mathcal A$ there exists a special decomposition of
$\theta$.
\end{theorem}

\begin{rem}
To characterize an efficient point, \textit{scalarization} methods are
commonly used. One of them consists of considering the following
optimization problem
\[
\sup\{\lambda (x) \dvtx  x\in\mathcal A\},
\]
where, denoting the topological dual of $Z$ by $Z^*$, $\lambda \in Z^*$
is such that $\lambda \geq0$ on $C$ and $\lambda >0$ on $C\setminus
\operatorname{lin}(C)$. If the optimum is attained, say at $\theta $,
then $\theta $ is $C$-efficient.

To see this, observe that if $x\in \operatorname{cone}(\mathcal A- \theta )\cap C$,
then we can write it as $x=k(w-\theta )$ for some $k\geq0$ and
$w\in\mathcal A$, and since $x\in C$ we see that either $k=0$ or
$w-\theta \in C$ in which case $\lambda (w)\geq\lambda (\theta )$ which
implies equality and that $w-\theta \in \operatorname{lin}(C)$.

More generally, for $\xi\in Z$ with $(\xi+C)\cap\mathcal A\neq
\varnothing$, the arg-max of the optimization problem
\[
\sup\{\lambda (x) ; x\in(\xi+C)\cap\mathcal A\},
\]
if it exists, is $C$-efficient.
\end{rem}

\begin{pf*}{Proof of Theorem~\ref{decomp}}
The proof uses scalarization.

Notice first that we only need to prove that there is a $\theta_0$ such
that
\[
\theta_0 \mbox{ is efficient in } K_0\cap(\theta-\mathcal
A_{1,T})\mbox{ with respect to }-\!\mathcal A_{1,T},
\]
with a general ${\mathcal F}_0$ (not necessarily trivial). This is
sufficient, since we may then apply the result to $\theta_{t,T}$ in an
inductive argument.

To make the scalarization argument we seek a linear function
\[
\lambda\dvtx S\rightarrow m{\mathcal F}_0,
\]
where $S=\operatorname{span}(K_0)=K_0-K_0$, with the properties that
%
\begin{equation}\label{scal}
\lambda\leq0\qquad\mbox{a.s. on } C\buildrel \mathrm{def}\over=\mathcal
A_{1,T}\cap S
\end{equation}
and
%
\begin{equation}\label{scal2}
[X\in C\mbox{ and } \lambda(X)=0\mbox{ a.s.}]\Rightarrow X\in
\operatorname{lin}(\mathcal A_{1,T}).
\end{equation}

First, notice that $\mathcal A_{1,T}$ is closed, and so is $S$ since it
is finitely ${\mathcal F}_0$-generated. Thus $C$ is a convex cone,
closed in $L^0({\mathcal F}_0;\mathbb R^d)$ and stable under
multiplication by elements of $b{\mathcal F}_0^+$. Moreover, since
$\mathcal A$ is arbitrage-free, so is $C$.

It follows from the abstract closed convex cone theorem of \cite{JBW}
that there exists a set-valued map $\Lambda\dvtx
\Omega\rightarrow{\mathcal P}(\mathbb R^d)$ such that:
\begin{longlist}[(3)]
\item[(1)] $\Lambda$ is almost surely a closed convex cone;

\item[(2)] $\Lambda$ is Effros--Borel measurable: that is, the event
$(\Lambda\cap U\neq\varnothing)$ is in ${\mathcal F}_0$ for any open
set $U\subset \mathbb R^d$;

\item[(3)] $C=\{X\in L^0({\mathcal F}_0,\mathbb R^d)\dvtx  X\in\Lambda
\mbox{ a.s}\}$.
\end{longlist}
It is easy to check that the map $\Lambda^*$, obtained by defining
$\Lambda^*(\omega)$ to be the polar cone of $\Lambda(\omega)$, also
satisfies (2) and is also almost surely a closed convex cone. It
follows from the fundamental measurability theorem of
\cite{Him} that there is a countable set $\{ Y_n\dvtx  n\geq1\}$ in
$L^0({\mathcal F}_0,\mathbb R^d)$ such that
%
\begin{equation}\label{scal4}
\Lambda^*(\omega)=\overline{\{ Y_n(\omega )\dvtx  n\geq1\}}.
\end{equation}

Now we claim that, setting
\[
\lambda=\sum2^{-n}Y_n/|Y_n|_{\mathbb R^d},
\]
$\lambda$ satisfies (\ref{scal}) and (\ref{scal2}).
Notice first that, since $\Lambda^*(\omega )$ is a.s. a closed convex cone,
%
\begin{equation}\label{scal3}
P(\lambda\in\Lambda^*)=1.
\end{equation}

To see that $\lambda$ satisfies (\ref{scal}): first take an $X\in C$,
then, by property (3), $X\in\Lambda\mbox{ a.s}$. Now, by (\ref{scal3}),
$X$ and $\lambda$ almost surely lie in polar cones in $\mathbb R^d$, so
(\ref{scal}) holds.

To prove that $\lambda$ satisfies (\ref{scal2}): suppose that $X\in C$
and $\lambda\cdot X=0$ a.s. It follows from the definition of $\lambda$
and (\ref{scal4}) that $Y_n\cdot X=0$ a.s., for each $n$. We can
conclude, again from (\ref{scal4}), that $\mu\cdot X=0$ a.s., for any
$\mu$ such that $\mu\in\Lambda^*$ a.s. Now this in turn implies that
$-X$ has the same property, which shows that $- X\in\Lambda$ a.s. We
conclude from (3) that $-X\in C$ and hence that $X\in
\operatorname{lin}(C)$.

Having obtained our linear function $\lambda$ which is negative on
$C\setminus \operatorname{lin}(C)$, we denote $K_0\cap(\theta-\mathcal
A_{1,T})$ by $\hat K_0$. Notice that $\hat K_0$ is closed since both
$K_0$ and $\mathcal A_{1,T}$ are. Now we claim that $\hat K_0$ is a.s.
bounded, that is, defining $M=\{|X|_{\mathbb R^d}\dvtx  X\in\hat
K_0\}$:
%
\begin{equation}\label{ess}
W^*\buildrel \mathrm{def}\over=\operatorname{ess} \operatorname{sup}
\{W\dvtx  W\in M\}<\infty\qquad\mbox{a.s.}
\end{equation}
To see this, notice first that $M$ is directed upward since, given $X$
and $Y$ in $\hat K_0$, $X\mathbf1_{(|X|_{\mathbb R^d}\geq|Y|_{\mathbb
R^d})}+Y\mathbf1_{(|X|_{\mathbb R^d}< |Y|_{\mathbb R^d})}\in\hat K_0$.
It follows that there is a sequence $(X_n)_{n\geq1}\subset\hat K_0$
such that
\[
|X_n|_{\mathbb R^d}\uparrow W^*\qquad\mbox{a.s.}
\]
Now define $F$ to be the event $(W^*=\infty)$.

Since $X_n\in\hat K_0$, there is a $Y_n\in\mathcal A_{1,T}$ with
$X_n+Y_n=\theta.$
Now, setting
\[
F_n=F\cap(|X_n|_{\mathbb R^d}>0)
\]
and multiplying by
$\frac{\mathbf1_{F_n}}{|X_n|_{\mathbb R^d}}$ we obtain:
\[
x_n+y_n=\tilde{\theta}\mathbf1_{F_n},
\]
where
\[
x_n=\mathbf1_{F_n}X_n/|X_n|_{\mathbb R^d},\qquad y_n=\mathbf
1_{F_n}Y_n/|X_n|_{\mathbb R^d}\quad\mbox{and}\quad
\tilde{\theta}=\theta /|X_n|_{\mathbb R^d}.
\]
Now, by Lemma A.2 in \cite{schacher}, we may take a strictly increasing
${\mathcal F}_0$-measurable random subsequence $(\tau_k)_{k\geq1}$ such
that $x_{\tau_k}$ converges almost surely, to $x$ say. From the
definition of $F_n$ we see that
$\tilde{\theta}\mathbf1_{F_n}\buildrel\mathrm{a.s.}\over
\longrightarrow 0$ and hence $y_n\buildrel\mathrm{a.s.}\over
\longrightarrow y$ for some $y\in\mathcal A_{1,T}$. But this implies
that $(x,y)$ is in $\mathcal N(K_0\times\mathcal A_{1,T})$ and hence
$x=0$ a.s. However, $|x|_{\mathbb R^d}=1$ a.s. on $F$ and so ${\mathbb
P}(F)=0$, establishing (\ref{ess}).

To complete the proof, observe that ${{\mathbb L}}\buildrel
\mathrm{def}\over=\{ \lambda\cdot X\dvtx  X\in\hat K_0\}$ is also
directed upward and so we may take a sequence
$(X_n)_{n\geq1}\subset\hat K_0$ such that
\[
\lambda\cdot X_n\uparrow l^* \buildrel
\mathrm{def}\over=\operatorname{ess}\operatorname{sup}{{\mathbb L}}.
\]
Now we take a strictly increasing ${\mathcal F}_0$-measurable random
sequence $(\sigma_k)_{k\geq1}$ such that $X_{\sigma_k}$ converges
almost surely, to $X$ say, and it follows from (\ref{ess}) that $X\in
L^0$.

Since $X_n\in\hat K_0$, there exists a $Y_n\in\mathcal A_{1,T}$ with
$\theta=X_n+Y_n$. Now
\[
Y_{\sigma_k}=\sum _{i=1}^\infty Y_i\mathbf1_{(\sigma_k=i)},
\]
and $\mathcal A_{1,T}$ is a closed convex cone, stable under
multiplication by elements of $m{\mathcal F}_0^+$, so
$Y_{\sigma_k}\in\mathcal A_{1,T}$. Similarly, $K_0$ is a closed convex
cone stable under multiplication by elements of $m{\mathcal F}_0^+$, so
it follows that $X_{\sigma_k}\in K_0 $ and, since
$X_{\sigma_k}+Y_{\sigma_k}=\theta$, that
\[
X_{\sigma_k}\in \hat K_0.
\]
By closure we deduce that $X\in\hat K_0$ and $\lambda\cdot X=l^*$ a.s.

Now the scalarization argument shows that we may take $\theta_0=X$
since if $Y\in\hat K_0$ and $X=Y+U$ with $U\in\mathcal A_{1,T}$ then
\[
\lambda\cdot X=\lambda\cdot Y+\lambda\cdot U,
\]
and since $U=X-Y\in\hat K_0-\hat K_0\subset \operatorname{span}(K_0)$
it follows that $U\in\mathcal A_{1,T}\cap S$ and so $\lambda\cdot
U\leq0$ a.s. But the maximality of $\lambda\cdot X$ now implies that
$\lambda\cdot U=0$ a.s., and we conclude from (\ref{scal}) that $U\in
\operatorname{lin}(\mathcal A_{1,T})$ which shows that $X$ is
efficient.
\end{pf*}

We shall now sketch a plan for the main result:
\begin{enumerate}[Step (2)]
\item[Step (1)] take a special decomposition $(\theta_0,\ldots,\theta
_T)$ for a maximal claim $\theta$;

\item[Step (2)] suppose that there exists a sequence $\mathbb
G=(G_1,\ldots,G_T)$ such that
%
\begin{equation}\label{adapt}
G_t\in{\mathcal F}_t\qquad\mbox{for each }t
\end{equation}
and
\begin{eqnarray}\label{null2}
& \mbox{whenever } y_t\in K_{t-1}-m{\mathcal F}_{t-1}^+\theta
_{t-1}\mbox{ with } -y_t(G_t^c)\in A_{t,T} &
\nonumber
\\[-8pt]
\\[-8pt]
\nonumber
&\mbox{ we can conclude that }y_t=0;&
\end{eqnarray}

\item[Step (3)] show that
%
\begin{equation}\label{proper}
\theta^{\mathbb G}\buildrel \mathrm{def}\over=
\theta_0+\theta_1(H_1)+\cdots+\theta_T(H_T),
\end{equation}
where $H_t\buildrel \mathrm{def}\over= G_1\cap\cdots\cap G_t$, is
properly maximal; to do this, show by backward induction that
%
\begin{equation}\label{proper2}
\theta^{\mathbb G}_{t,T}\mbox{ is properly maximal in }A_{t,T};
\end{equation}

\item[Step (4)] Show that, using randomization, there exists a
sequence $(\mathbb G^n)_{n\geq1}$ such
that each $\mathbb G^n$ satisfies properties (\ref{adapt}) and (\ref
{null2}) and ${\mathbb P}(G_t^n)\uparrow1$
for each $t$.
\end{enumerate}
For the rest of this section we assume that
$\operatorname{lin}(\mathcal A)=\{0\}$.

We now implement Step (3).
For the initial step in the induction we need the following result:

\begin{lem}\label{l0}
Suppose that $K$ is a finitely ${\mathcal F}$-generated convex cone and
is arbitrage-free.

Let $\xi\in K$. Then,
%
\begin{equation}\label{clos}
\overline{[K,\xi]}=K-m{\mathcal F}^+ \xi
\end{equation}
and hence is finitely generated. Moreover, \textup{[}using the ordering
cone $L^{0,+}_d({\mathcal F})$\textup{]}
\begin{equation}\label{cmax}
\xi\mbox{ is maximal in $K$ if and only if it is properly maximal.}
\end{equation}
\end{lem}

\begin{pf}
Suppose that $\lambda\in m{\mathcal F}^+$ and define $\lambda_n=\min
(\lambda,n)$, then $(n-\lambda_n)\xi\in K$ and hence $-\lambda_n
\xi=(n-\lambda_n)\xi- n\xi\in[K,\xi]$. Hence $K-m{\mathcal F}^+
\xi\subset \overline{[K,\xi]}$. Conversely, since $K-m{\mathcal F}^+
\xi$ is finitely generated it is closed and contains $[K,\xi]$,
so~(\ref{clos}) is satisfied.

To prove (\ref{cmax}), suppose $\xi$ is maximal in $K$ and that $y$ is
an arbitrage in $\overline{[K,\xi]}$, so that $y=x-\alpha \xi$ with
$x\in K$ and $\alpha \in m{\mathcal F}^+$. It follows that $\widebar
{x}\buildrel \mathrm{def}\over=\frac{1}{\alpha } x\mathbf1_{(\alpha
>0)}\in K$ and $\widebar{x}=\frac{1}{\alpha } y\mathbf1_{(\alpha
>0)}+\xi\mathbf1_{(\alpha >0)}$. Hence, since $y\geq0$,
\[
z\buildrel \mathrm{def}\over=\widebar{x}+\xi\mathbf1_{(\alpha =0)} \geq
\xi
\]
and $z\in K$. Since $\xi$ is maximal we get
$\widebar{x}+\xi\mathbf1_{(\alpha =0)}=z=\xi$ and then $y\mathbf
1_{(\alpha >0)}=0$. Finally, since $\mathbf1_{(\alpha
=0)}y=\mathbf1_{(\alpha =0)}x\in K$ and $K$ is arbitrage-free, we
conclude that
\[
\mathbf1_{(\alpha =0)}y=0
\]
and hence that $y=0$.
\end{pf}

\begin{theorem}\label{proper1}
Suppose that $\operatorname{lin}(\mathcal A)=\{0\}$, that
$\theta\in\mathcal A$ is maximal, that
$\theta=\theta_0+\cdots+\theta_T$ is a special decomposition of
$\theta$, that $\mathbb G$ satisfies \textup{(\ref{adapt})} and
\textup{(\ref{null2})} and that $\theta^{\mathbb G}$ is as defined in
\textup{(\ref{proper})}. Then
\[
\theta^{\mathbb G}\mbox{ is properly maximal in }\mathcal A.
\]
\end{theorem}

\begin{pf}
As announced, we shall show that (\ref{proper2}) holds for each $t$.

Assume that $\theta^{\mathbb G}_{t+1,T}\mbox{ is properly maximal in
}\mathcal A_{t+1,T}$. Now it is easy to check that
\[
\overline{[\mathcal A_{t,T},\theta^{\mathbb G}_{t,T}]}=
\overline{\overline{[K_t,\theta^{\mathbb G}_t]}+\overline{[\mathcal
A_{t+1,T},\theta^{\mathbb G}_{t+1,T}]}},
\]
so if we can show that
\[
S_t\buildrel \mathrm{def}\over=\overline{[K_t,\theta^{\mathbb
G}_t]}+\overline{[\mathcal A_{t+1,T},\theta^{\mathbb G}_{t+1,T}]}
\]
is closed and arbitrage-free then the inductive step is complete. Then
Lemma \ref{l0} gives us the initial step (for $S_T$).

($S_t$ \textit{is closed})

We do this by showing that
\[
N\buildrel \mathrm{def}\over=\mathcal N(\overline{[K_t,\theta^{\mathbb
G}_t]}\times\overline{[\mathcal A_{t+1,T},\theta^{\mathbb
G}_{t+1,T}]})=\{\underline{0}\},
\]
and appealing to Lemma \ref{s}. To do this, notice first that
(\ref{clos}) tells us that $\overline{[K_t,\theta^{\mathbb
G}_t]}=K_t-m{\mathcal F}_{t}^+\theta ^{\mathbb G}_{t}\subset K_t-
m{\mathcal F}_{t}^+\theta_{t}$. Now notice that if $z\in
\overline{[\mathcal A_{t+1,T},\theta^{\mathbb G}_{t+1,T}]}$ then,
taking a sequence $z_n\in[\mathcal A_{t+1,T},\theta^{\mathbb
G}_{t+1,T}]$ converging to $z$ we see that, since $G_{t+1}^c\subset
H_{t+1}^c$ and $\theta^{\mathbb G}_{t+1,T}$ is supported on $H_{t+1}$,
$z_n\mathbf1_{G_{t+1}^c}\in\mathcal A_{t+1,T}$ for each $n$, and hence
$z\mathbf1_{G_{t+1}^c}\in\mathcal A_{t+1,T}$.

So if $(y,z)\in N$ then
\[
y(G_{t+1}^c)+z(G_{t+1}^c)=0
\]
and so it follows from (\ref{null2}) that $y=0$ and hence that $z=0$.

($S_t$ \textit{is arbitrage-free})

Suppose that $f$ is an arbitrage in $S_t$, that is, $f\geq0$ and
$f=y+z$ with $y\in\overline{[K_t,\theta^{\mathbb G}_t]}$ and
$z\in\overline{[\mathcal A_{t+1,T},\theta^{\mathbb G}_{t+1,T}]}$. Then
\[
0=y+(z-f)
\]
and so (since $L^{0,-}_d\subset\mathcal A_{t+1,T}$), $(y,z-f)$ is in $N$
and so $y=z-f=0$. It follows that $z=f\geq0$ and since, by the
inductive hypothesis,
$\overline{[\mathcal A_{t+1,T},\theta^{\mathbb G}_{t+1,T}]}=S_{t+1}$ is
arbitrage-free, the inductive step follows.
\end{pf}

We implement the Step (4) of the proof plan as follows:

First define
\[
\hat\Omega _i=\mathbb N,\qquad \hat\sigma_t=2^{\Omega _t},\qquad
\hat\Omega =\hat\Omega _1\times\cdots\times\hat\Omega _T,\qquad
\hat{\mathcal F}_t=\sigma_1\otimes\cdots\otimes\sigma_t,
\]
and then define
\[
\tilde\Omega =\Omega \times\hat\Omega\quad\mbox{and}\quad \tilde
{\mathcal F} _t={\mathcal F}_t\otimes\hat{\mathcal F}_t.
\]
To complete the randomization, define a probability measure $\tilde
{\mathbb P} $ on $(\tilde\Omega ,\tilde{\mathcal F} _T)$ by setting
\[
\tilde{\mathbb P} ={\mathbb P}\otimes\hat{\mathbb P} \otimes\cdots
\otimes\hat{\mathbb P} ,
\]
where $\hat{\mathbb P} $ is the probability measure on $\mathbb N$
defined by $\hat{\mathbb P} (\{k\})=2^{-k}$.

Now set
%
\begin{equation}\label{rand}
G^n_t=\Omega \times\hat\Omega _{1}\times\cdots\times \hat\Omega
_{t-1}\times\{1,\ldots,n\}\times\hat\Omega _{t+1}\times\cdots\times
\hat\Omega _{T}.
\end{equation}
It is clear that $G^n_t\uparrow\tilde\Omega $ as $n\uparrow\infty$
for each $t$.

We extend then the definition of the cone $\mathcal A$ to the new
setting by defining $\tilde K_t$ to be the convex cone
${\tilde{\mathcal F}}_t$-generated by the same generators as $K_t$,
that is, if
\[
K_t= \Biggl\{\sum_{i=1}^n \alpha_i \Pi^i_t\dvtx  \alpha_i\in m{\mathcal
F}_t^+ \Biggr\}
\]
then
\[
\tilde K _t= \Biggl\{\sum_{i=1}^n \alpha_i \Pi^i_t\dvtx  \alpha_i\in
m\tilde {\mathcal F} _t^+ \Biggr\};
\]
and then set $\tilde{\mathcal A}=\tilde K_0+\cdots+\tilde K_T$.

\begin{lem}\label{randlem}
Under Assumption \textup{\ref{neat}}, the convex cone $\tilde{\mathcal
A}$ is closed, arbitrage-free and the null strategies subset $\mathcal
N(\tilde K_0\times\cdots\times\tilde K_T)$ is trivial. Moreover each
maximal claim in $\mathcal A$ is also maximal in $\tilde{\mathcal A}$.
\end{lem}

\begin{pf}
Each property of $\tilde{\mathcal A} $ follows from the corresponding
property for $\mathcal A$ in the same way. So, for example, the null
strategies for $\mathcal A$ form a vector space, $N$~say. Now take
$(\xi_0,\ldots,\xi_T)\in\tilde N$, where $\tilde N$ is the collection
of null strategies for $\tilde\mathcal K$, then fix
$(i_1,\ldots,i_T)\in\mathbb N^T$ then
$(\xi_0(\cdot),\xi_1(\cdot;i_1),\ldots,\xi_T(\cdot,i_1,\allowbreak \ldots ,\allowbreak i_T))\in
N$ and so
\[
-(\xi_0(\cdot),\xi_1(\cdot;i_1),\ldots,\xi_T(\cdot,i_1,\ldots
,i_T))\in N,
\]
and since $(i_1,\ldots,i_T)$ is arbitrary, $-(\xi_0,\ldots,\xi
_T)\in\tilde N$ and hence $\tilde N$ is a vector space. The same
method---of freezing those arguments of an $\tilde{\mathcal F}
_t$-measurable random variable which are in $\hat\Omega $ will
establish each of the results.
\end{pf}

We need one more lemma before we can give the main result:

\begin{lem}\label{null4}
If $\theta=\theta_0+\cdots+\theta_T$ is a special decomposition of
$\theta\in\mathcal A$, then, for each $t$, the null strategies subset
$\mathcal N(\overline{[K_t,\theta _t]}\times\mathcal A_{t+1,T})$ is
trivial.
\end{lem}

\begin{pf}
Since $\theta=\theta_0+\cdots+\theta_T$ is a special decomposition of
$\theta$, it follows that, defining
\[
\hat K_t=K_t\cap(\theta_{t,T}-\mathcal A_{t+1,T}),
\]
$\theta_t$ is efficient in $\hat K_t$ with respect to $-\mathcal
A_{t+1,T}$, that is,
%
\begin{equation}\label{null6}
(\hat K_t-\theta_t)\cap(-\mathcal A_{t+1,T})=\{0\}.
\end{equation}

Now we know from Lemma \ref{l0} that $\overline{[K_t,\theta_t]}=K_t-
m{\mathcal F}_t^+\theta_t$, so any null strategy for
$\overline{[K_t,\theta_t]}\times A_{t+1,T}$ is of the form
$(x_t-\lambda_t\theta_t,x_{t+1,T})$, where $x_t\in K_t$, $\lambda
_t\in m{\mathcal F}_t^+$ and
$x_{t+1,T}\in A_{t+1,T}$. Now, take such a triple, so that
%
\begin{equation}\label{null3}
x_t-\lambda_t\theta_t+x_{t+1,T}=0,
\end{equation}
and multiply (\ref{null3}) by $\mathbf1_{(\lambda_t=0)}$ to get:
\[
x_t\mathbf1_{(\lambda_t=0)}+x_{t+1,T}\mathbf1_{(\lambda_t=0)}=0.
\]
So, we conclude that
%
\begin{equation}\label{null7}
x_t\mathbf1_{(\lambda_t=0)}=0,
\end{equation}
because
$\mathcal N(\mathcal K)=\{\underline0\}$ and so $\mathcal N(K_t\times
A_{t+1,T})=\{\underline0\}$.

Now multiply (\ref{null3}) by $\alpha_t\buildrel
\mathrm{def}\over=\frac{1}{\lambda_t}\mathbf1_{(\lambda _t>0)}$ to
obtain
%
\begin{equation}\label{null5}
\hspace*{5mm} \alpha_tx_t+\theta_t\mathbf1_{(\lambda_t=0)}-\theta_t
+\alpha_t x_{t+1,T}=\alpha_tx_t-\theta_t\mathbf1_{(\lambda
_t>0)}+\alpha_t x_{t+1,T}=0.
\end{equation}
Now $\alpha_t\in m{\mathcal F}_t^+$ so $\alpha_tx_t\in K_t$ and, since
$\theta_t\in K_t$, we see that
\[
y_t\buildrel
\mathrm{def}\over=\alpha_tx_t+\theta_t\mathbf1_{(\lambda _t=0)}\in K_t.
\]
Moreover, from (\ref{null5})
\[
y_t=\theta_t-\alpha_tx_{t+1,T}=\theta_{t,T}-
(\alpha_tx_{t+1,T}+\theta_{t+1,T}),
\]
and so
\[
y_t\in(\theta_{t,T}-\mathcal A_{t+1,T}).
\]
We deduce that
\[
y_t\in\hat K_t.
\]
Now $y_t-\theta_t\in\hat K_t -\theta_t$ and $y_t-\theta_t\in- \mathcal
A_{t+1,T}$ so we deduce from (\ref{null6}) that
\[
y_t-\theta_t=0\mbox{ which implies that } \alpha_tx_t-\theta
_t\mathbf1_{(\lambda_t>0)}=0,
\]
and, multiplying by $\lambda_t$ and adding (\ref{null7})
we obtain the desired result that
\[
x_t-\lambda_t\theta_t=0.
\]\upqed
\end{pf}

\begin{theorem}\label{t1manu}
Let $\theta \in\mathcal A$ be a maximal claim in $\mathcal A$ (or
indeed in $\tilde{\mathcal A}$). Then there exists a sequence of
properly maximal claims $(\theta ^n)_{n\geq 1}$ in $\tilde{\mathcal A}$
which converge a.s. to $\theta $.
\end{theorem}

\begin{pf}
Thanks to Lemma \ref{randlem} we may work with $\tilde{\mathcal A}$
thoughout. We fix the special decomposition
$\theta=\theta_0+\cdots+\theta_T$ and, taking $\mathbb G^n$ as in
(\ref{rand}), define $\theta^n=\theta^{\mathbb G^n}$
using~(\ref{proper}).

Now suppose that
\[
y\in\overline{[\tilde K _{t-1},\theta_{t-1}]}\quad\mbox{and}\quad
{-}y\mathbf1_{(G^n_t)^c}\in\tilde{\mathcal A} _{t,T},
\]
then, setting $z=-y$,
\[
z\mathbf1_{(G^n_t)^c}\in\tilde{\mathcal A} _{t,T}
\]
and
\[
y\mathbf1_{(G^n_t)^c}+z\mathbf1_{(G^n_t)^c}=0.
\]
Now take any $j>n$ then, since $y$ is $\tilde{\mathcal F} _{t-1}$-measurable
and $z$ is $\tilde{\mathcal F} _t$-measurable,
\[
y(\omega ;\hat\omega _1,\ldots,\hat\omega _{t-1})+z(\omega ;\hat \omega
_1,\ldots,\hat\omega _{t- 1},j)=0\qquad\mbox{a.s.}
\]
Finally, taking $(\hat\omega _1,\ldots,\hat\omega
_{t-1})=(i_1,\ldots,i_{t-1})$ we see that
\[
y(\cdot,i_1,\ldots,i_{t-1})\in K_{t-1}-m{\mathcal F}_{t-1}^+\theta
_{t-1}=\overline{[K_{t-
1},\theta_{t-1}]}
\]
and
\[
z(\cdot,i_1,\ldots,i_{t-1},j)\in\mathcal A_{t,T}
\]
for each choice of $i_1,\ldots,i_{t-1},j$ and so it follows from Lemma
\ref{null4} that $y(\cdot,i_1,\ldots,i_{t-1})=0$ for each choice of
$i_1,\ldots,i_{t-1},j$ and so $y=0$. The fact that $\theta^n$ is
properly maximal now follows from Theorem \ref{proper1}. It is obvious
that $\theta^n\buildrel\mathrm{a.s.}\over\longrightarrow \theta$ as
$n\uparrow\infty$.
\end{pf}

\begin{rem}
Since the convergence in Theorem \ref{t1manu} follows from a truncation,
it is clear that if the special decomposition used has the property that
$\theta_t\in L^p({\mathcal F}_t,\mathbb R^d)$ for each $t$ then
convergence of the properly maximal sequence will also be
in $L^p$ by the dominated convergence theorem.
\end{rem}

\section{The case $\operatorname{lin}(\mathcal A)\neq\{0\}$}

In the case where $\operatorname{lin}(\mathcal A)\neq\{0\}$ we may
still assume that $\mathcal N(\mathcal K)=\{\underline{0}\}$, however
the conclusion of Lemma \ref{null4} fails, that is, we may no longer
conclude that, with $\theta_t$ being the $t$th component of a special
decomposition of $\theta$, $\mathcal N(\overline{[K_t,\theta_t]}\times
\mathcal A_{t+1,T})=\{\underline{0}\}$.

The way around this problem is to focus on $t=0$ and define $\sim$, an
equivalence relation on elements of $K_0\cap(\theta-\mathcal A_{1,T})$,
as follows:
\[
x\sim y\Leftrightarrow x-y\in \operatorname{lin}(\mathcal A_{1,T}).
\]

\begin{rem}
Notice that if $\theta_0$ is efficient in $K_0\cap(\theta-\mathcal
A_{1,T})$ (with respect to $-\mathcal A_{1,T}$), then every element of
the equivalence class $[\theta_0]$ is efficient. To see this, take
$z\in [\theta_0]$, so $z\in K_0\cap(\theta-\mathcal A_{1,T})$ and
$z-\theta_0\in- \mathcal A_{1,T}$.
\end{rem}

Now we can easily show that, defining
\[
\Sigma_t=m{\mathcal F}_t^+[\theta_t],
\]
the correct generalization of Lemma \ref{null4} holds.

\begin{lem}\label{null11}
If $\theta_t$ is efficient, the null space $\mathcal N((K_t-\Sigma
_t)\times \mathcal A_{t+1,T} )$ is a vector space.
\end{lem}

\begin{pf}
As indicated, we need only to prove the result in the case where
$t=0$, provided we do not assume that ${\mathcal F}_0$ is trivial.

Suppose $x\in K_0$, $\xi\sim\theta_0$, $\lambda\in m{\mathcal
F}_0^+$, $z\in
\mathcal A_{1,T}$ and
%
\begin{equation}
x-\lambda\xi+z=0.
\end{equation}
It is immediate that
%
\begin{equation}\label{null8}
(x+\xi)-(1+\lambda) \xi+z=0,
\end{equation}
and, dividing (\ref{null8}) by $1+\lambda$ we get
\[
\tilde x-\xi+\tilde z=0.
\]
It follows, since $\tilde x\in K_0\cap(\theta -\mathcal A_{1,T})$ and
$\xi$ is efficient, that $\tilde z\in \operatorname{lin}(\mathcal
A_{1,T})$ and therefore that $\tilde x\sim\xi$. And so $\tilde
x\sim\theta_0$ and $z\in \operatorname{lin}(\mathcal A_{1,T})$. So
$\xi-\tilde x\in K_0-\Sigma_0$ and, multiplying by $1+\lambda$,
\[
\lambda\xi-x\in K_0-\Sigma_0
\]
and so $\mathcal N((K_0-\Sigma_0)\times\mathcal A_{1,T} )$ is a
vector space.
\end{pf}

We now have another problem since Lemma \ref{l0} is no longer apparently
relevant---at first sight it does not look as though $K_0-\Sigma_0$ is finitely
generated, so it is not clear that it is closed.

\begin{lem}\label{fin2}
For each $t$, there is a $\xi_t\in[\theta_t]$ such that
%
\begin{equation}\label{spec}
K_t-\Sigma_t=\overline{[K_t,\xi_t]}=K_t-m{\mathcal F}_t^+\xi_t
\end{equation}
and so $K_t-\Sigma_t$ is closed.
\end{lem}

\begin{pf}
As before, we only need to prove the lemma for $t=0$ and a
nontrivial~${\mathcal F}_0$. Now, since $\Sigma_0\subset K_0$, it is
clear that
%
\begin{equation}\label{max}
K_0-\Sigma_0=K_0+\Sigma_0-\Sigma_0.
\end{equation}
We shall prove that, for the right choice of $\xi_0\in\Sigma_0$,
%
\begin{equation}\label{max2}
\Sigma_0-\Sigma_0=(K_0-m{\mathcal F}_0^+\xi_0)\cap\bigl(m{\mathcal
F}_0\xi_0+\operatorname{lin}(A_{1,T})\bigr)=\Sigma_0-m{\mathcal
F}_0^+\xi_0,
\end{equation}
by
showing that
%
\begin{equation}\label{max3}
\Sigma_0-\Sigma_0\subset(K_0-m{\mathcal F}_0^+\xi_0)\cap
\bigl(m{\mathcal
F}_0\xi_0+\operatorname{lin}(A_{1,T})\bigr)\subset\Sigma_0- m{\mathcal
F}_0^+\xi_0.
\end{equation}
Notice that if (\ref{max3}) holds then there must be equality
throughout, since $m{\mathcal F}_0^+\xi_0\subset\Sigma_0$, and the
result will then follow immediately from (\ref{max}) and (\ref{max2}).

We define $\xi_0$ as follows.

First recall that the generators of $K_0$ are $(\Pi^i_0)_{1\leq
i\leq m}$. Now define
\[
\Phi\buildrel \mathrm{def}\over=\Biggl\{(\alpha_1,\ldots,\alpha_m)\dvtx
\sum _i\alpha_i\Pi^i_0\in\Sigma_0; \alpha_i\in m{\mathcal F}_0^+,
|\alpha_i|\leq1 \mbox{ for }i=1,\ldots,m \Biggr\}.
\]
It is clear that $\Phi$ is a convex set, closed in $L^0({\mathcal
F}_0;\mathbb R^m)$.

Now define $p\dvtx \Phi\rightarrow\mathbb R^+$ by
\[
p(\underline\alpha)=\sum _{i=1}^m{\mathbb P}(\alpha_i>0).
\]
Denote $\sup _{\alpha\in\Phi}p(\alpha)$ by $p^*$ (notice that $p^*\leq
m$) and take a sequence $(\underline{\alpha_n})_{n\geq 1}\subset\Phi$
such that $p(\underline{\alpha_n})\uparrow p^*$. It follows from the
convexity and closure of $\Phi$ that
\[
\sum_{k=1}^n
2^{-k}\underline{\alpha_k}+2^{-n}\underline{\alpha_{n+1}}\buildrel
\mathrm{a.s.}\over\longrightarrow \sum_{k=1}^\infty
2^{-k}\underline{\alpha_k}\buildrel \mathrm{def}\over=\underline{\hat
\alpha }\in\Phi
\]
and
\[
p(\hat\alpha )=p^*.
\]
Now define
\[
\xi_0\buildrel \mathrm{def}\over=\sum_{i=1}^m \hat\alpha _i\Pi_0^i.
\]
The convexity and closure of $\Phi$ ensures that $\xi_0\in\Sigma_0$.
Notice that it follows from the definition of $\hat\alpha $ that if
$x=\sum_{i=1}^m \alpha_i\Pi^i_0\in \Sigma_0$ then
\[
{\mathbb P}\bigl((\alpha_i>0)\cap(\hat\alpha
_i=0)\bigr)=0\qquad\mbox{for each }i.
\]

Denote the middle term in (\ref{max2}) by $R$.

($\Sigma_0-\Sigma_0\subset R$)

Since $\Sigma_0\subset K_0$ and $\psi\sim\phi\Rightarrow\psi-\phi \in
\operatorname{lin}(\mathcal A_{1,T})$, which implies that
\[
\Sigma_0\subset\bigl(m{\mathcal
F}_0\xi_0+\operatorname{lin}(A_{1,T})\bigr),
\]
we see that
\[
\Sigma_0\subset R.
\]
Now take $x\in\Sigma_0$, so $x\in K_0$ and $z\buildrel
\mathrm{def}\over=x- \alpha^+\xi_0\in \operatorname{lin}(A_{1,T})$ for
some $\alpha^+\in m{\mathcal F}_0^+$. It follows that
$-x=-z-\alpha^+\xi_0$ so $-x\in(m{\mathcal F}_0\xi
_0+\operatorname{lin}(A_{1,T}))$. All that remains for this step is to
prove that
%
\begin{equation}\label{max8}
-x\in(K_0-m{\mathcal F}_0^+\xi_0).
\end{equation}
Recall that $\xi_0$ has maximal support in $\Sigma_0$, so if $x=\sum _i
\alpha_i \Pi_0^i$ and $\xi_0=\sum _i\hat\alpha _i \Pi _0^i$ then
$\beta\buildrel \mathrm{def}\over= \max _i\{\frac{\alpha_i}{\hat\alpha
_i}\}<\infty$ a.s. Since $\beta\in m{\mathcal F}_0^+$ it follows that
\[
\beta\xi_0-x\in K_0
\]
and hence, expressing $-x$ as
$(\beta\xi_0-x)-\beta\xi_0$, we conclude that (\ref{max8})
holds.\vadjust{\goodbreak}

($R\subset\Sigma_0-m{\mathcal F}_0^+\xi_0$)

Take $y\in R$. Since $y\in(K_0-m{\mathcal F}_0^+\xi_0)$ we may write
it as
\[
y=x-\alpha^+\xi_0,
\]
with $x\in K_0$ and $\alpha^+\in m{\mathcal F}_0^+$. Moreover, since
$y\in (m{\mathcal F}_0\xi_0+\operatorname{lin}(\mathcal A_{1,T}))$ we
may write it as
\[
y=\gamma\xi_0+z,
\]
with $\gamma\in m{\mathcal F}_0$ and $z\in \operatorname{lin}(\mathcal
A_{1,T})$. Denoting the positive and negative parts of $\gamma$ by
$\gamma^+$ and $\gamma^-$ respectively, it follows that
%
\begin{equation}\label{max4}
x+\gamma^-\xi_0=(\gamma^++\alpha^+)\xi_0+z.
\end{equation}
The right-hand side of (\ref{max4}) is clearly in $(m{\mathcal
F}_0^+\xi_0+\operatorname{lin}(\mathcal A_{1,T}))$ and the left-hand
side is clearly in $K_0$, so we conclude that the common value, $w$
say, is in $\Sigma_0$.

Finally, observe that
\[
y=w-(\alpha^++\gamma^-)\xi_0,
\]
and so $y\in\Sigma_0-m{\mathcal F}_0^+\xi_0$.
\end{pf}

Now we suitably generalize condition (\ref{null2}) and
Theorem~\ref{proper1}.

Suppose that $\theta\in\mathcal A$ and it is decomposed as
$\theta=\theta_0+\cdots+\theta_T$.

\begin{defi}
For each $t$, define $N_t$ as the projection onto the first component
of the nullspace $\mathcal N((K_t-\Sigma_t)\times\mathcal A_{t+1,T})$
and define $N_t^\perp$ as the orthogonal complement of $N_t$ (this is
well defined thanks to Lemma~\textup{\ref{null11}} and Lemma
\textup{A.4} in~\textup{\cite{schacher})}.
\end{defi}

Further suppose that there exists a sequence $\mathbb
G=(G_1,\ldots,G_T)$ such that
%
\begin{equation}\label{adapt2}
G_t\in{\mathcal F}_t\qquad\mbox{for each }t
\end{equation}
and
%
\begin{eqnarray}\label{null9}
&\mbox{whenever } y_t\in(K_{t-1}-m{\mathcal F}_{t-1}^+[\theta
_{t-1}])\cap N_{t-1}^{\perp}& \nonumber
\\[-8pt]
\\[-8pt]
\nonumber &\mbox{ with } -y_t(G_t^c)\in\mathcal A_{t,T}\mbox{ we may
conclude that }y_t=0.&
\end{eqnarray}

\begin{theorem}\label{proper3}
Suppose that $\theta\in\mathcal A$ is maximal, that
$\theta=\xi_0+\cdots+\xi_T$ is a special decomposition of $\theta$ with
$\underline\xi$ as in Lemma \textup{\ref{fin2}}, so that
\[
K_{t}-m{\mathcal F}_{t}^+[\xi_{t}]=K_{t}-m{\mathcal F}_{t}^+\xi_{t}.
\]
Suppose, in addition, that $\mathbb G$ satisfies
\textup{(\ref{adapt2})} and \textup{(\ref{null9})} and that
$\theta^{\mathbb G}$ is as defined in~\textup{(\ref{proper})}, then
\[
\theta^{\mathbb G}\mbox{ is properly maximal in }\mathcal A.
\]
\end{theorem}

\begin{pf}
The argument mirrors the proof of Theorem \ref{proper1}.

As before we need to show that
\[
S_t\buildrel \mathrm{def}\over=\overline{[K_t,\xi_t(H_t)]}\cap
N_t^\perp+\overline{[\mathcal A_{t+1,T},\theta^{\mathbb G}_{t+1,T}]}
\]
is closed and arbitrage-free.

($S_t$ \textit{is closed})

We do this by showing that
\[
N\buildrel \mathrm{def}\over=\mathcal N\bigl(\bigl({[K_t,\xi_t(H_t)]}\cap
N_t^\perp \bigr)\times\overline{[\mathcal A_{t+1,T},\theta^{\mathbb
G}_{t+1,T}]} \,\bigr)=\{\underline{0}\}.
\]
To do this, notice first that (\ref{clos}) and (\ref{spec}) tell us that
\[
\overline{[K_t,\xi_t(H_t)]}=K_t-m{\mathcal F}_{t}^+[\xi _t(H_t)]\subset
K_t- m{\mathcal F}_{t}^+\xi_{t}.
\]
Now notice that, as before, if $z\in \overline{[\mathcal
A_{t+1,T},\theta^{\mathbb G}_{t+1,T}]}$ then,
$z\mathbf1_{G_{t+1}^c}\in\mathcal A_{t+1,T}$.

So if $(y,z)\in N$ then
\[
y(G_{t+1}^c)+z(G_{t+1}^c)=0
\]
and so it follows from (\ref{null9}) that $y=0$ and hence that $z=0$.

($S_t$ \textit{is arbitrage-free})

The argument is unchanged.
\end{pf}

The proof of the revised version of Theorem \ref{t1manu} is essentially
unchanged. Since the statement does not involve
$\operatorname{lin}(\mathcal A)$ we do not repeat it.

\section{Further comments}

A slight modification of Theorem \ref{proper3} states, under some mild
assumptions, that for any maximal claim $\theta\in\mathcal A$, there exists
a sequence of properly maximal claims $\theta_n$ which converges to
$\theta$ in probability.

\begin{theorem}\label{proper4}
Given $\theta\in\mathcal A$ is maximal, take a special decomposition of
$\theta\dvtx  \theta=\xi_0+\cdots+\xi_T$, with $\underline\xi$ as in
Lemma~\textup{\ref{fin2}}, so that
\[
K_{t}-m{\mathcal F}_{t}^+[\xi_{t}]=K_{t}-m{\mathcal F}_{t}^+\xi_{t}.
\]
Suppose there exists a sequence $\mathbb G^n$ satisfying
\textup{(\ref{adapt2})} and \textup{(\ref{null9})}, with each $G^n_t$
converging to~$\Omega$.

Now define the sequence $\theta^n\buildrel
\mathrm{def}\over=\theta^{\mathbb G^n}$ as in \textup{(\ref{proper})}.
Then
\[
\mbox{ the sequence }\theta^n\mbox{ is properly maximal in }\mathcal A
\]
and
\[
\theta_n\rightarrow\theta\qquad\mbox{in probability}.
\]
\end{theorem}

Unfortunately we are unable to construct such a sequence $\mathbb G^n$
in a general setting. We have adopted a randomization approach that
allows us to construct such sequence.

We remark that hedging such a randomized sequence is still possible
``in the market without randomization.'' By this we mean that, since
trades in the randomized market take place at the same bid-ask prices
as in the original market, an individual trader may perform the
randomizations and hedge accordingly in the original market.

\section*{Acknowledgments}
The authors are grateful for many fruitful discussions with Jon Warren
on the topics of this paper. We also thank two anonymous referees for
some very helpful suggestions.

\printaddresses

\end{document}